\documentclass[12pt]{article}
\usepackage{amsmath, amssymb, amsthm}
\usepackage{mathrsfs}
\usepackage{amssymb}
\usepackage{amsmath}
\usepackage{hyperref}
\usepackage{marvosym}
\usepackage{pict2e}
\usepackage{graphicx}
\usepackage{geometry}
\usepackage{color}
\def\q{\hfill\rule{1ex}{1ex}}
\def\0{\emptyset}

\def\p{{\bf Proof.}\quad}
\def\q{\hfill\rule{1ex}{1ex}}

\newtheorem{lem}{Lemma}

\newtheorem{theo}{Theorem}[section]
\newtheorem{coro}[theo]{Corollary}
\newtheorem{conj}[theo]{Conjecture}

\begin{document}

\title{\bf Chords of longest cycles in graphs with large circumferences}
\author{Haidong Wu~~~~~Shunzhe Zhang\\}

\author{Haidong Wu\footnote{E-mail: hwu@olemiss.edu}
~~~~~Shunzhe Zhang\footnote{E-mail: szhang12@go.olemiss.edu}\\
Department of Mathematics\\
The University of Mississippi\\
University, MS 38677, U.S.A\\
}

\date{}
\maketitle
\baselineskip 16.7pt \setcounter{page}{1}

\begin{abstract}

A long-standing conjecture of Thomassen says that every longest cycle of a $3$-connected graph has a chord. Thomassen (2018) proved that if $G$ is a $2$-connected cubic graph, then any longest cycle must have a chord. He also showed that in any 3-connected graph with minimum degree at least four, some longest cycle must contain a chord. Harvey proved that every longest cycle has a chord for graphs with a large minimum degree. He also conjectured that any longest cycle in a 2-connected graph with minimum degree at least three has a chord. In this paper, we prove that both Thomassen's and Harvey's conjectures are true for graphs with large circumferences. We also prove a more general result for the existence of chords in longest cycles containing a linear forest. 
\vskip 0.1cm

\noindent {\bf Keywords:} longest cycles, chords, $3$-connected, planar graph

\end{abstract}

\section{Introduction}

In 1976, Thomassen raised the following famous conjecture:

\begin{conj}\label{conj1.1} {\rm (Thomassen~\cite{Th1989}, also in 
\cite[(Conjecture 8.1)]{Al1985})}
\label{originalconj}
Every longest cycle of a $3$-connected graph has a chord.
\end{conj}

The general conjecture remains open, although many partial results have been proven (see, \cite{Bi2008}, \cite{Ha2017}, \cite{Ka2007}, \cite{Li2024}, \cite{LiZh20031}, \cite{LiZh20032}, \cite{Th1989}, \cite{Th1997}, \cite{Th2018}, \cite{WY2025}, \cite{Wu-Zhang} \cite{Zh1987}). Zhang~\cite{Zh1987} (1987) showed that the conjecture holds for 3-connected planar graphs that are either cubic or have minimum degree of at least four. 

\begin{theo}\label{theo1.2} {\rm (Zhang~\cite{Zh1987})}
Let $G$ be a $3$-connected planar graph which is either cubic or with minimum degree at least $4$. Then any longest cycle of $G$ must have a chord. 
\end{theo}

Thomassen~\cite{Th1997} showed that Conjecture \ref{originalconj} holds when $G$ is 3-connected cubic. In 2018, he further extended this result \cite{Th2018}.

\begin{theo}\label{theo1.3} {\rm (Thomassen~\cite{Th2018})}
Every longest cycle in a $2$-connected cubic graph has a chord. Moreover, for any $3$-connected graph with minimum degree at least four,  some longest cycle must have a chord. 
\end{theo}

Harvey~\cite{Ha2017} proved that every longest cycle has a chord for graphs with large minimum degree. 
He also made a more general conjecture than Thomassen's conjecture, and the following is a special case of that general conjecture. 

\begin{conj}\label{harvey} {\rm (Harvey~\cite{Ha2017})}
Let $G$ be a $2$-connected graph such that the minimum degree is at least three. Then every longest cycle has a chord. 
\end{conj}

 Wang and Yue \cite{WY2025} showed that Thomassen's conjecture is true for $3$-connected graphs with circumference at least $n-5$ where $n$ is the number of vertices of the graph. In this paper, we show that both Thomassen's and Harvey's conjectures are true for graphs with large circumferences. Note that we do not need to assume the connectivity condition. 

\begin{theo}\label{main1} 
Let $G$ be a simple graph with $n$ vertices and having minimum degree of at least three. Suppose that the circumference of $G$ is at least $n-\frac{1+\sqrt{4n-3}}{2}$. Then every longest cycle of $G$ has a chord. 
\end{theo}

By Theorem \ref{main1} and the fact $\sqrt{n}\leq \frac{1+\sqrt{4n-3}}{2}$, we obtain the following corollary.

\begin{coro}
Let $G$ be a simple graph with $n$ vertices and having minimum degree of at least three. Suppose that the circumference of $G$ is at least $n-\sqrt{n}$. Then every longest cycle of $G$ has a chord. 
\end{coro}

The lower bound in our theorem \ref{main1} can likely be improved.  Harvey gave an example in \cite[Page 2]{Ha2017} showing a class of connected graphs with minimum degree $\sqrt{n}-1$, but with a longest cycle having size $\sqrt{n}$ (here $n=|V(G)|$) having no chord. Similarly, we can construct a graph as follows (See Figure 1). Let $C_{5k}$ be a cycle with a partition set $X_1,\ldots,X_k$ of $V(C_{5k})$, where $X_i=\{x_{5i-4},x_{5i-3},x_{5i-2},x_{5i-1},x_{5i}\}$ and $k\ge 1$. For each $i\in \{1,\ldots,k\}$, three vertices $x_{5i-4},x_{5i-2},x_{5i}$ are adjacent to $y_i$ which is not in $C_{5k}$ and each $x_{5i-3},x_{5i-1}$ is identified with a vertex of a copy of $K_4$. Then the resulting graph has a longest cycle having no chord with the circumference being $c(G)={5n\over 12}$. Another example is to start with a wheel with $k+1$ vertices. Subdivide each rim edge into two edges first, then for each vertex of degree $2$, take a copy of $N=K_4$ and identify a vertex of $N$ with the vertex of degree two.  These examples show that the best possible lower bound for $c(G)$ one can hope for to ensure that any longest cycle containing a chord is linear.

\begin{figure}[h]
\centering
\includegraphics[width=0.9\textwidth,height=0.48\textwidth]{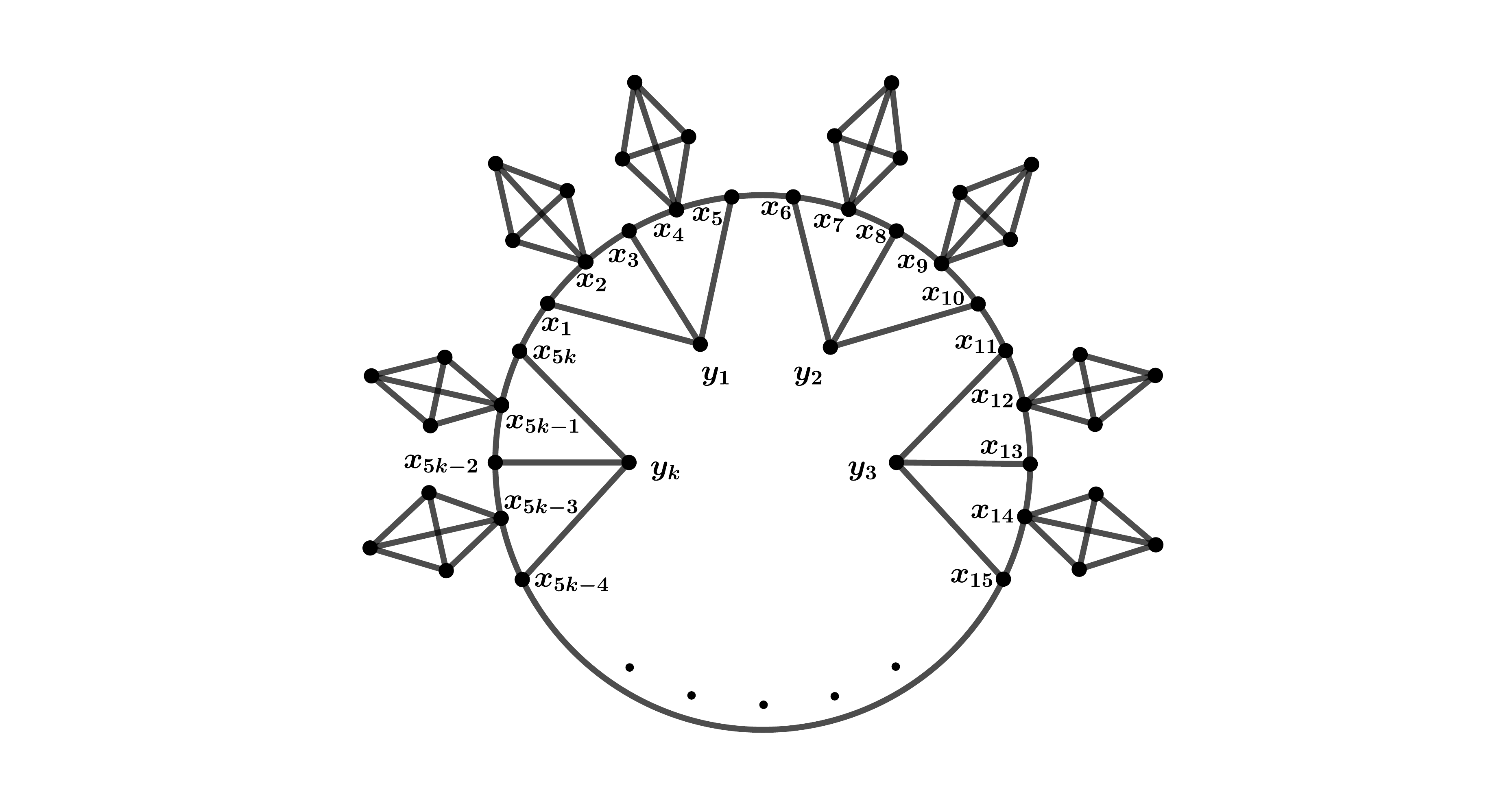} 
\vspace{-.5cm}
\caption{A graph having a longest cycle with no chord and with circumference $5n\over 12$.}      
\label{fig1}   
\end{figure}

\noindent {\bf Question 1}:  Let $G$ be a simple graph with $n$ vertices and with minimum degree at least three. Find the best possible constant $c$  such that if the circumference of $G$ is at least $cn$, then every longest cycle in $G$ contains a chord.

Our examples described above show that $c>{5\over 12}$. A linear forest is a graph such that each component is a path. A trivial linear forest is the empty set. In \cite{Wu-Zhang}, we made the following conjecture:

\begin{conj}\label{Wu} 
Let $G$ be a $k$-connected graph ($k\ge 2$) and let $F$ be a linear forest subgraph of $G$ with $l$ edges and $t$ isolated vertices such that $l+t\le k-2$. Then every longest cycle of $G$ passing through $F$ has a chord.
\end{conj}

Note that the connectivity condition in Conjecture \ref{Wu} ensures the existence of a cycle containing $F$. For the special case of $F$ consisting of a single edge,  we consider the longest cycles (if they exist) containing a specified edge and show that Conjecture \ref{Wu} is true for connected graphs with large circumferences.

\begin{theo}\label{main2} 
Let $G$ be a simple graph with $n$ vertices and having minimum degree of at least three, and $e$ be an edge of $G$. Suppose that the length of the longest cycles containing $e$ is at least $n-\frac{1+\sqrt{4n-3}}{2}+1$. Then every longest cycle of $G$ containing $e$ has a chord. 
\end{theo}

The following corollary follows immediately from Theorem \ref{main2}. 

\begin{coro}
Let $G$ be a simple graph with $n$ vertices and having minimum degree of at least three, and $e$ be a specified edge of $G$. Suppose that the longest cycle containing $e$ has length at least  $n-\sqrt{n}+1$.  Then every longest cycle of $G$ containing $e$ has a chord. 
\end{coro}

Furthermore, we consider the longest cycles (if they exist) containing a linear forest with at most one edge (that is, $F$ consists of isolated vertices plus possibly one edge)  and show that Conjecture \ref{Wu} is true when the circumference is large.

\begin{theo}\label{main3} 
Let $F$ be a linear forest with at most one edge of a simple graph $G$ with $n$ vertices having minimum degree of at least three. Suppose that the length of the longest cycle containing $F$ is at least $n-\sqrt{n-1}+1$. Then every longest cycle of $G$ containing $F$  has a chord. 
\end{theo}

In \cite{Ha2017}, Harvey also made the following conjecture:

\begin{conj} {\rm (Harvey~\cite{Ha2017})}
If $\delta(G)\ge \sqrt{n}$ (where $n$ is the number of vertices of $G$). Then every cycle of maximum order in $G$ contains a chord. 
\end{conj}

Both Thomassen's conjecture and Harvey's conjecture \ref{harvey} are very hard. It would be interesting to know if the lower bound in Theorem \ref{main1} can be improved to $2\sqrt{n}$ if we assume the graph is $2$-connected. This is weaker than Harvey's conjecture (Conjecture \ref{harvey}).

\noindent {\bf Question 2}: Let $G$ be a $2$-connected graph with  $n$ vertices and with minimum degree at least three such that the circumference of $G$ is at least $2\sqrt{n}$. Does every longest cycle of $G$ contain a chord?

Proofs of our main results will be given in Section $2$. We use Bondy and Murty~\cite{Bo2008} for terminology and notation not defined here. Let $C$ be a cycle with an arbitrary orientation. If $u,v\in V(C)$, we use $C[u,v]$ (or $uCv$) and $\overline{C}[u,v]$ (or $u\overline{C}v$) to denote the subpath of $C$ from $u$ to $v$ along the orientation of $C$ and the same subpath in reverse order, respectively. Set $C(u,v]=C[u,v]\setminus \{u\}$, $C[u,v)=C[u,v]\setminus \{v\}$ and $C(u,v)=C[u,v]\setminus \{u,v\}$.

\section{Proofs}

In this section, we give proofs of our theorems. First, we list the following two elementary lemmas, which will be used in our proofs.

\begin{lem}\label{lemma1}
Let $C=v_0v_1\ldots v_tv_0$ be a cycle of $G$ containing a linear forest $F$ with at most one possible edge $e$, and $u, v$ be two vertices outside of $C$.  If $u$ is adjacent to two vertices $v_i$ and $v_j$ ($1\le i+1<j\le t$) and $v$ is adjacent to both $v_{i+1}$ and $v_{j+1}$, where the index is read modular $t$ such that $e$ is neither $v_iv_{i+1}$ nor $v_jv_{j+1}$.  Then $G$ has a longer cycle than $C$, and this cycle contains $F$.   
\end{lem}

\p The new cycle $v_iuv_jv_{j-1}\ldots v_{i+1}vv_{j+1}\ldots v_{i-1}v_i$ is a cycle with length $|V(C)|+2$ and apparently it contains all vertices of $C$ and thus contains all vertices of $F$. Moreover, as $e$ is neither $v_iv_{i+1}$ nor $v_jv_{j+1}$, this new cycle contains $e$ also. \q

\begin{lem}\label{lemma2}
Let $n\geq 3t-a\geq 1$, where $n,t,a$ are positive integers. Then $n-\frac{1+\sqrt{4n-3}}{2}>2t-a-1$. Furthermore, if $t\geq 4$ and $a\geq 2$, then $n-\frac{1+\sqrt{4n-3}}{2}>2t-a$.
\end{lem}

\p Let $f(n)=n-\frac{1+\sqrt{4n-3}}{2}$. It is easily checked that $f(n)$ is increasing for $n\geq 1$. Since $n\geq 3t-a\geq 1$, we deduce that $f(n)\geq f(3t-a)=3t-a-\frac{1+\sqrt{12t-(4a+3)}}{2}=2t-a-1+\frac{2(t-1)^2+2a}{2t+1+\sqrt{12t-(4a+3)}}>2t-a-1$. Furthermore, if $t\geq 4$ and $a\geq 2$, then $f(n)\geq f(3t-a)=3t-a-\frac{1+\sqrt{12t-(4a+3)}}{2}=2t-a+\frac{2(t-1)(t-3)+2(a-2)}{2t-1+\sqrt{12t-(4a+3)}}>2t-a$.
\q

Next we prove our first main result, Theorem \ref{main1}. Our first theorem will follow from the next theorem.

\begin{theo}\label{main1-conn}
Let $G$ be a simple connected graph with $n$ vertices and having minimum degree of at least three. Suppose that the circumference of $G$ is at least $n-\frac{1+\sqrt{4n-3}}{2}$. Then every longest cycle of $G$ has a chord. 
\end{theo}

{\bf Proof}. Let $G$ be a simple connected graph with a longest cycle $C$ having size $n-k$ where $n=|V(G)|$ and $k\le \frac{1+\sqrt{4n-3}}{2}$ be a positive integer. Then $n\ge k^2-k+1$. Let $f(n)=n-\frac{1+\sqrt{4n-3}}{2}$. By our assumption, $|V(C)|\geq f(n)$. On the contrary, suppose that $C$ has no chord. Assume that $G-V(C)$ has $t$ components. We orient $C$ first, then contract each component of $G-V(C)$ into a single vertex and remove multiple edges. Then the new graph $G_1$ is connected with $t\le k$ vertices not in $V(C)$, and $C$ is still a chordless longest cycle of $G_1$. Moreover, each vertex of $C$ still has degree at least three. Apparently, $t\ge 1$. If $t=1$, then as each vertex of the cycle $C$ has degree at least three, $G_1$ must be a wheel. This is a contradiction as $C$ is not a longest cycle of $G_1$. Thus $2\le t\le k$. Next, we show that $G_1$ must have $ n-k+t \geq t^2-t+1$ vertices. Indeed, if $n-k+t\le t^2-t$, then $k^2-k+1\le n\le t^2-2t+k$, thus $k^2-2k< t^2-2t$. As $t, k\ge 2$, we conclude that $k<t$, a contradiction. Let $v$ be a vertex in $S=V(G_1-V(C))$ with maximum degree in $G_1$. Color all neighbors of $v$ with red color and all immediate neighbors of these red vertices along the orientation of $C$ as blue (these vertices are not adjacent to $v$, apparently). Suppose that $v$ has at least $t$ red neighbors in $C$. By the Pigeonhole principle, there is at least one vertex in $S-\{v\}$ (which has $t-1$ vertices) adjacent to two blue vertices of $C$ (as each vertex of $C$ has degree at least three in $G_1$). This is a contradiction by Lemma \ref{lemma1}. Therefore, the degree of $v$ in $G_1$ is at most $t-1$. It follows that that $|V(C)|\leq t(t-1)$. Thus $t\geq 3$ as $|V(C)|\geq 3$.

Now, we will prove that the degree of $v$ in $G_1$ is at most $t-2$. On the contrary, assume that the degree of $v$ in $G_1$ is $t-1$. By a similar argument to that in the previous paragraph, each vertex in $S-\{v\}$ is adjacent to exactly one blue vertex of $C$. Let $X=\{x_0,\ldots,x_{t-2}\}$ and $Y=\{y_0,\ldots,y_{t-2}\}$ denote the sets of all red and blue vertices of $C$, respectively, where $x_iy_i\in E(C)$ for $i\in \{0,\ldots,t-2\}$. Let $S-\{v\}=\{v_0,\ldots,v_{t-2}\}$ such that $y_i\in N_{G_1}(v_i)$ for $i\in \{0,\ldots,t-2\}$. Apparently, $|V(C)|\geq 2t-2$. Next, we will show that $|V(C)|=2t-2$. Otherwise, $|V(C[x_i,x_{i+1}))|\geq 3$ for some $i\in \{0,\ldots,t-2\}$, where the index is read modular $t-1$. Without loss of generality, assume that $|V(C[x_0,x_1))|\geq 3$. Let $z_0$ be the immediate neighbor of $y_0$ along the orientation of $C$. Since each vertex of $C$ has degree at least three in $G_1$, $z_0$ is adjacent to $v_i$ for some $i\in \{1,\ldots,t-2\}$. But then $x_0vx_i\overline{C}z_0v_iy_iCx_0$ is a cycle with length $|V(C)|+1$ in $G_1$, a contradiction.

Clearly, $n=|V(G)|\geq |V(G_1)|\geq |V(C)|+t=3t-2$. If $t\geq 4$, then by Lemma \ref{lemma2}, $f(n)>2t-2=|V(C)|$, which contradicts to the fact that $|V(C)|\geq f(n)$. If $t=3$, then $|V(C)|=4$, $|S|=3$ and thus $|V(G_1)|=7$. So $n-\frac{1+\sqrt{4n-3}}{2}=f(n)\leq |V(C)|=4$. It follows that $3\leq n\leq 7$. This implies that $|V(G)|=|V(G_1)|=7$ and $G\cong G_1$. But then $d_G(v)=d_{G_1}(v)=t-1=2$, a contradiction.

Therefore, the degree of $v$ in $G_1$ is at most $t-2$, and thus $G_1$ has at most $t(t-2)+t=t^2-t$ vertices, a contradiction to the fact that $|V(G_1)|\geq t^2-t+1$. This contradiction shows that $C$ must have a chord.  \q

{\bf Proof of \ref{main1}}: By Theorem \ref{main1-conn}, we need only prove the theorem for disconnected graphs. Let $G$ be a simple disconnected graph with $n$ vertices and having minimum degree of at least three. Suppose that the circumference of $G$ is at least $n-\frac{1+\sqrt{4n-3}}{2}$, and $C$ be a longest cycle of $G$. Let $H$ be the connected component of $G$ containing $C$ and suppose that $|V(H)|=p$. Then $n>p$. As $n-\frac{1+\sqrt{4n-3}}{2}$ is a strictly increasing function for $n\ge 2$, we deduce that $n-\frac{1+\sqrt{4n-3}}{2}>p-\frac{1+\sqrt{4p-3}}{2}$. Thus, the circumference of $H$ is larger than 
$p-\frac{1+\sqrt{4p-3}}{2}$. By Theorem \ref{main1-conn}, the longest cycle $C$ of $H$ has a chord, thus $C$ has a chord in $G$. \q

Next we prove our second main result, Theorem \ref{main2}. We prove the following result for connected graphs first.

\begin{theo}\label{main2-conn}
Let $G$ be a simple connected graph with $n$ vertices and having minimum degree of at least three, and $e$ be an edge of $G$. Suppose that the length of the longest cycle containing $e$ is at least $n-\frac{1+\sqrt{4n-3}}{2}+1$.  Then every longest cycle of $G$ containing $e$ has a chord. 
\end{theo}

{\bf Proof}. Let $C$ be a longest cycle containing an edge $e$ (which must exist) of $G$ having size $n-k$ where $n=|V(G)|$ and $k<\frac{1+\sqrt{4n-3}}{2}$ be a positive integer. Then $n>k^2-k+1$. Let $g(n)=n-\frac{1+\sqrt{4n-3}}{2}+1$. By our assumption, $|V(C)|\geq g(n)$. On the contrary, suppose that $C$ has no chord. Assume that $G-V(C)$ has $t$ components. We orient $C$ first, then contract each component of $G-V(C)$ into a single vertex and remove the multiple edges. Then the new graph $G_1$ is connected with $t\le k$ vertices not in $V(C)$, and $C$ is still a chordless longest cycle of $G_1$ containing $e$. Moreover, each vertex of $C$ still has degree at least three. If $t=1$, then as each vertex of the cycle $C$ has degree at least three, $G_1$ must be a wheel. This is a contradiction as $C$ is not a longest cycle of $G_1$ containing $e$. Thus $2\le t\le k$. Next, we show that $G_1$ must have $n-k+t> t^2-t+1$ vertices. Indeed, if $n-k+t\le t^2-t+1$, then $k^2-k+1<n\le t^2-2t+k+1$, thus $k^2-2k< t^2-2t$. As $t, k\ge 2$, we conclude that $k<t$, a contradiction.

Suppose that $e=xy$ is an edge of $C$, where $x, y\in V(C)$. First we assume that for any $v\in S=V(G_1-V(C))$, $v$ is not adjacent to both $x$ and $y$. Color all neighbors of $v$ with red colors. Reorient $C$ if necessary, we can assume that none of the immediate vertices of the red vertices along the orientation of $C$ are in the set $\{x,y\}$. Color these immediate neighbors of the red vertices along the orientation of $C$ as blue. Thus, these blue vertices are not adjacent to $v$. Suppose that $v$ has at least $t$ red neighbors in $C$. By the Pigeonhole principle, there is at least one vertex in $S-\{v\}$ (which has $t-1$ vertices) adjacent to two blue vertices of $C$ (as each vertex of $C$ has degree at least three in $G_1$). This is a contradiction by Lemma \ref{lemma1}. Therefore, the degree of $v$ in $G_1$ is at most $t-1$. It follows that $|V(C)|\leq t(t-1)$. Thus $t\geq 3$ as $|V(C)|\geq 3$.

\noindent {\bf Claim 1}  The degree of $v$ in $G_1$ is at most $t-2$. 

On the contrary, assume that the degree of $v$ in $G_1$ is $t-1$. Then $v$ has $t-1$ red neighbors and there are $t-1$ blue vertices in $C$. Since each blue vertex has degree at least three in $G_1$, and $C$ is chordless, each blue vertex is connected to some vertex of $S-\{v\}$. Using Lemma \ref{lemma1}, we deduce that each vertex in $S-\{v\}$ is adjacent to exactly one blue vertex of $C$.  Let $X=\{x_0,\ldots,x_{t-2}\}$ and $Y=\{y_0,\ldots,y_{t-2}\}$ denote the sets of all red and blue vertices of $C$, respectively, where $x_iy_i\in E(C)$ for $i\in \{0,\ldots,t-2\}$. Let $S-\{v\}=\{v_0,\ldots,v_{t-2}\}$ such that $y_i\in N_{G_1}(v_i)$ for $i\in \{0,\ldots,t-2\}$. Without loss of generality, assume that $e\in E(x_0Cx_1)$. Apparently, $|V(C)|\geq 2t-2$. Next, we will show that $2t-2\leq |V(C)|\leq 2t-1$. Suppose that $|V(C)|\geq 2t$. Then either there is a segment, say $C[x_k,x_{k+1})$ having at least four vertices, or there are at least two segments $C[x_i,x_{i+1})$, $C[x_j,x_{j+1})$, $(i\neq j)$ each having exactly three vertices, where the index is read modular $t$. In the former case, we consider the case when $k=0$. For other $k\in\{1, 2, \ldots, t-2\}$, the argument still works (indeed, it is even simpler since $C[x_k,x_{k+1})$ does not contain $e$ then).  So assume that $|V(C[x_0,x_1))|\geq 4$.  Then reorient $C$ if necessary, we can assume that the immediate vertex $z_0$ of $y_0$ along the orientation of $C$ satisfies $e\notin \{x_0y_0,y_0z_0\}$. Since $d_{G_1}(z_0)\geq 3$ and $C$ is a longest cycle containing $e$, $z_0$ is adjacent to $v_i$ for some $i\in \{1,\ldots,t-2\}$. But then $x_0vx_i\overline{C}z_0v_iy_iCx_0$ is a cycle containing $e$ with length $|V(C)|+1$ in $G_1$, a contradiction. Now, in the latter case, we assume that there are at least two segments $C[x_i,x_{i+1})$, $C[x_j,x_{j+1})$, $(i\neq j)$ each having exactly three vertices. Therefore, at least one of them, say $C[x_i,x_{i+1})$ does not contain $e$. Let $z_i$ be the immediate neighbor of $y_i$ along the orientation of $C$. Since $d_{G_1}(z_i)\geq 3$ and $C$ is a longest cycle containing $e$, $z_i$ is adjacent to $v_j$ for some $j\in \{0,\ldots,t-2\}\setminus \{i\}$. But then $x_ivx_j\overline{C}z_iv_jy_jCx_i$ is a cycle containing $e$ with length $|V(C)|+1$ in $G_1$, a contradiction. Hence, $2t-2\leq |V(C)|\leq 2t-1$. Let $|V(C)|=2t-a$, where $a\in \{1,2\}$. Clearly, $n=|V(G)|\geq |V(G_1)|\geq |V(C)|+t=3t-a$. By Lemma \ref{lemma2}, we deduce that $g(n)=(n-\frac{1+\sqrt{4n-3}}{2})+1>2t-a=|V(C)|$, which contradicts the fact that $|V(C)|\geq g(n)$. This completes the proof of Claim 1. \q

By Claim 1, the degree of $v$ in $G_1$ is at most $t-2$. As $v$ is an arbitrary vertex in $S$, we deduce that $G_1$ has at most $t(t-2)+t=t^2-t$ vertices, a contradiction to the fact that $|V(G_1)|>t^2-t+1$. This contradiction shows that $C$ must have a chord.

Now suppose that there is at least one vertex $u\in S=V(G_1-V(C))$ which is adjacent to both $x$ and $y$. Color all neighbors of $u$ with red colors and assume that $y$ is the immediate neighbor of $x$ along the orientation of $C$. Color these immediate neighbors of the red vertices along the orientation of $C$ as blue, except the immediate neighbor $y$ of $x$. Thus, these blue vertices are not adjacent to $u$. Suppose that $u$ has at least $t+1$ red neighbors in $C$. Then there are at least $t$ blue vertices in $C$. By the Pigeonhole principle, there is at least one vertex in $S-\{u\}$ (which has $t-1$ vertices) adjacent to two blue vertices of $C$ (as each vertex of $C$ has degree at least three in $G_1$). This is a contradiction as we can get a longer cycle containing $e$ by Lemma \ref{lemma1}. Therefore, the degree of $u$ in $G_1$ is at most $t$. If $t=2$, then $|S|=2$. Let $S=\{u,v\}$. Since $d_{G_1}(u)\leq 2$ and $u$ is adjacent to both $x$ and $y$, we have $N_{G_1}(u)=\{x,y\}$. But then $|V(C)|=3$, as otherwise, $v$ must be adjacent to all vertices of $C-\{x, y\}$, and thus $G_1$ has a cycle containing $e$ with length $|V(C)|+1$ in $G_1$, a contradiction. We conclude that $|V(G_1)|=5$ and $n-\frac{1+\sqrt{4n-3}}{2}+1=g(n)\leq |V(C)|=3$. It follows that $3-\sqrt{2}\leq n\leq 3+\sqrt{2}$. But then $|V(G)|=n<5=|V(G_1)|$, a contradiction. Hence, $t\geq 3$.

\noindent {\bf Claim 2}  The degree of $u$ in $G_1$ is at most $t-1$. 

By the last paragraph, the degree of $u$ in $G_1$ is at most $t$. So if Claim 2 is not true, then the degree of $u$ in $G_1$ is exactly equal to $t$.  Therefore $u$ has exactly $t$ red neighbors in $C$, two of which are $x$ and $y$ and $e=xy$,  and just like the proof in Claim 1,   there are exactly $t-1$ blue vertices in $C$. 

By a similar argument to that before, each vertex in $S-\{u\}$ is adjacent to exactly one blue vertex of $C$. Let $X=\{x_0,\ldots,x_{t-1}\}$ and $Y=\{y_1,\ldots,y_{t-1}\}$ denote the sets of all red and blue vertices of $C$, respectively, where $x_iy_i\in E(C)$ for $i\in \{1,\ldots,t-1\}$ and $x_0=x$, and $x_1=y$. Let $S-\{u\}=\{u_1,\ldots,u_{t-1}\}$ such that $y_i\in N_{G_1}(u_i)$ for $i\in \{1,\ldots,t-1\}$. Thus $e\in E[x_0Cx_1]=xy$, where $x_1=y$. Apparently, $|V(C)|\geq 2t-1$. Next, we will show that $|V(C)|=2t-1$. Suppose that $|V(C)|\geq 2t$. Then there is a segment, say $C[x_k,x_{k+1})$ for some $k$ ($1\le k\le t-1$)  having at least three vertices, where the index is read modular $t$. Then $e\notin C[x_k,x_{k+1})$.  Let $z_k$ be the immediate neighbor of $y_k$ in $C$.  Since $d_{G_1}(z_k)\geq 3$ and $C$ is a longest cycle containing $e$, $z_k$ is adjacent to $u_i$ for some $i\in \{1,\ldots,t-1\}$ different from $k$. But then $x_kux_i\overline{C}z_ku_iy_iCx_k$ is a cycle containing $e$ with length $|V(C)|+1$ in $G_1$, a contradiction.  Hence, $|V(C)|=2t-1$. Clearly, $n=|V(G)|\geq |V(G_1)|\geq |V(C)|+t=3t-1$. By Lemma \ref{lemma2}, we deduce that $g(n)=(n-\frac{1+\sqrt{4n-3}}{2})+1>2t-1=|V(C)|$, which contradicts that $|V(C)|\geq g(n)$. Therefore, the degree of $u$ in $G_1$ is at most $t-1$. This completes the proof of Claim 2. \q

Now for any $v\in S-\{u\}$, by a similar argument to those in Claims 1 and 2, we can show that $v$ is adjacent to at most $t-2$ vertices of $V(C)-\{x, y\}$ in $G_1$. Thus $G_1$ has at most $t-1+(t-1)(t-2)+t=t^2-t+1$ vertices, a contradiction to the fact that $|V(G_1)|>t^2-t+1$. This contradiction shows that $C$ must have a chord. \q

 {\bf Proof of \ref{main2}}: By Theorem \ref{main2-conn}, we need only prove the theorem for disconnected graphs. Let $G$ be a simple disconnected graph with $n$ vertices and having minimum degree of at least three, and $e$ be a specified edge. Suppose that the size of longest cycles of $G$ containing $e$ is at least $n-\frac{1+\sqrt{4n-3}}{2}+1$, and $C$ be a longest cycle of $G$ containing $e$. Let $H$ be the connected component of $G$ containing $C$ and suppose that $|V(H)|=p$. Then $n>p$. As $n-\frac{1+\sqrt{4n-3}}{2}+1$ is a strictly increasing function for $n\ge 2$, we deduce that $n-\frac{1+\sqrt{4n-3}}{2}+1>p-\frac{1+\sqrt{4p-3}}{2}+1$. Thus, the size of the longest cycles in $H$ containing $e$  is larger than $p-\frac{1+\sqrt{4p-3}}{2}+1$. By Theorem  \ref{main2-conn}, the longest cycle $C$ of $H$ containing $e$ has a chord, thus $C$ has a chord in $G$. \q

Next we prove our third main result, Theorem \ref{main3}. Our result follows from the next result,  and we will omit the straightforward proof for disconnected graphs.

\begin{theo}
Let $F$ be a linear forest with at most one edge of a simple connected graph $G$ with $n$ vertices having minimum degree of at least three. Suppose that the length of the longest cycle containing $F$ is at least $n-\sqrt{n-1}+1$. Then every longest cycle of $G$ containing $F$  has a chord. 
\end{theo}

{\bf Proof}. Let $F$ be a linear forest of $G$. We only need to show that the theorem is true when there exists a cycle containing $F$.  Suppose that $C$ is a longest cycle containing $F$ having size $n-k$ where $k<\sqrt{n-1}$ is a positive integer. Then $n> k^2+1$. Suppose that $C$ has no chord. We orient $C$ first, then contract each component of $G-V(C)$ into a single vertex and remove the multiple edges. Then the new graph $G_1$ is connected with $t\le k$  vertices (where $t$ is the number of components of  $G-V(C)$ ) not in $V(C)$, and $C$ is still a chordless longest cycle of $G_1$ containing $F$. Moreover, each vertex of $C$ still has degree at least three. If $t=1$, then as each vertex of the cycle $C$ has degree at least three, $G_1$ must be a wheel. This is a contradiction as $C$ is not a longest cycle of $G_1$ containing $F$. Thus $2\le t\le k$. Next, we show that $G_1$ must have $n-k+t> t^2+1$ vertices. Indeed, if $n-k+t\le t^2+1$, then $k^2+1<n\le t^2-t+k+1$, thus $k^2-k< t^2-t$. As $t, k\ge 1$, we conclude that $k<t$, a contradiction.   

Suppose that $e=xy$, where $x, y\in V(C)$. First we assume that any vertex $v$ of $S=V(G_1-V(C))$ is not adjacent to both $x$ and $y$. Color all neighbors of $v$ with red colors. Reorient $C$ if necessary, we can assume that none of the immediate vertices of the red vertices along the orientation of $C$ are in the set $\{x, y\}$. Color these immediate neighbors of the red vertices along the orientation of $C$ as blue.  Thus, these blue vertices are not adjacent to $v$. Suppose that $v$ has at least $t$ red neighbors in $C$. By the Pigeonhole principle, there is at least one vertex in $S-\{v\}$ (which has $t-1$ vertices) adjacent to two blue vertices of $C$ (as each vertex of $C$ has degree at least three in $G_1$). By Lemma \ref{lemma1}, $G_1$ has a longer cycle containing $F$, a contradiction. Therefore, the degree of $v$ in $G_1$ is at most $t-1$, and thus $G_1$ has at most $t(t-1)+t=t^2$ vertices, a contradiction.

Now suppose that there is a vertex $u\in S=V(G_1-V(C))$ adjacent to both $x$ and $y$. Using a similar argument to that in the previous paragraph, we can show that (1) the degree of $u$ in $G_1$ is at most $t$, and (2)  any $v\in S-\{u\}$ is adjacent to at most $t-1$ vertices of $V(C)-\{x, y\}$. Thus $G_1$ has at most $t+(t-1)(t-1)+t=t^2+1$ vertices. This is a contradiction to the fact that $G_1$ has more than $t^2+1$ vertices. Therefore $C$ must have a chord.  \q

Let $C$ be a circuit of a matroid and $e$ be a non-loop element not in $C$. We say that $e$ is a {\it chord} of $C$ if $e$ is in the closure of $C$. A natural question is whether our main results can be extended to cographic matroids. Let $G$ be a simple connected graph with $m$ edges and $n$ vertices. It is well known that any largest bond of $G$ has size at most $p(G)=m-n+2$. If $B$ is a bond of $G$ with size of at least $p(G)-1=m-n+1$, then at least one of the two components of $G-B$ is a tree, thus $B$ has a chord for $M^*(G)$. However, the following example shows that the conclusion may not be true if $|B|\le p(G)-2$. Take two disjoint cycles of length $n$ with vertex sets $X$ and $Y$ respectively. For each vertex $v$ of $X$, add an edge from $v$ to each vertex of $Y$ to form a graph $G$. Let $B$ be the bond induced by the partition $(X, Y)$. Then it is easily checked that  $B$ is a largest bond of $G$ such that $|B|=n^2=p(G)-2$, and $B$ has no chord in the cographic matroid $M^*(G)$.

\section{\bf Acknowledgments} 

{ H. Wu's research is supported by the 2025 CLA grant of the University of Mississippi. }

\end{document}